\documentclass[oneside]{amsart}
\usepackage[T1]{fontenc}
\usepackage[utf8]{inputenc}
\usepackage{amsthm}

\makeatletter
\numberwithin{equation}{section}
\numberwithin{figure}{section}
\theoremstyle{plain}
\newtheorem{thm}{\protect\theoremname}

\usepackage{algorithmic}

\makeatother

\providecommand{\theoremname}{Theorem}

\begin{document}

\title{Counterexamples to Theorem 1 of Turyn's and Storer's Paper ``On
Binary Sequences'' }

\author{Jürgen Willms}

\email{willms.juergen@fh-swf.de}

\address{Institut für Computer Science, Vision and Computational Intelligence,
Fachhochschule Südwestfalen, D-59872 Meschede, Germany}
\begin{abstract}
Explicit counterexamples to Theorem 1 of R. Turyn's and J. Storer's
often-cited paper “On Binary Sequences” \cite{turyn1961binary} are
given. Theorem 2 of their paper is well known; it states that Barker
sequences of odd length $n>13$ do not exist. Since the proof of Theorem
2 relies on the incorrect Theorem 1, the proof of the often-cited
Theorem 2 as presented in \cite{turyn1961binary} is therefore not
correct. It is not at all clear how Theorem 1 and its proof must be
modified in \cite{turyn1961binary} in order to get a correct proof
of the well-known Theorem 2.
\end{abstract}



\date{7.5.2014 }

\maketitle

\section{Introduction}

Barker sequences are a prominent example of a family of binary sequences
which have specific aperiodic autocorrelation properties. There is
also a wide range of engineering applications where Barker sequences
are used. For example, both the IEEE 802.11b wireless LAN standard
and the GPS satellite navigation system use a Barker sequence for
modulation purposes. A Barker sequence $x=(x_{1},x_{2},\cdots,x_{n})$
is a binary sequence with the property that each of its aperiodic
autocorrelation $c_{k}=\sum_{i=1}^{n-k}x_{i}x_{i+k}$ is in magnitude
as small as possible, i. e. is either 0, 1 or -1. Theorem 2 of Turyn's
and J. Storer's well-known paper “On Binary Sequences” \cite{turyn1961binary}
states that there exists no Barker sequence for odd sequence length
$n>13.$ According to \cite{borwein2008barker} ``their proof is
elementary, though somewhat complicated''. The proof of Theorem 2
is based on Theorem 1 which is of a more technical nature. In the
following, we will give explicit counterexamples to Turyn's and Storer's
Theorem 1 \emph{(iv)}. As a consequence, the proof of the well-known
Theorem 2 in \cite{turyn1961binary} is not correct.

\section{Theorem 1 of Turyn's and Storer's Paper ``On Binary Sequences''}

In the following $x$ will as in \cite{turyn1961binary} always denote
a binary sequence with length $n$. Thus we have $x=(x_{1},x_{2},\cdots,x_{n})$
with $x_{i}\in\{-1,1\}$ for all $i=1,\cdots,n$. As in \cite{turyn1961binary}
we will say that a binary sequence $x=(x_{1},x_{2},\cdots,x_{n})$
will satisfy equation \emph{(k) }for $1\leq k<\frac{n-1}{2}$ if 
\[
(k)\qquad\qquad\qquad\qquad\frac{1+(-1)^{k+1}}{2}=\sum_{i=1}^{k}(-1)^{i+1}x_{i}x_{2k+2-i}.
\]

Note that whether the sequence $x$ satisfies equation \emph{(k)}
depends only on the first $2k+1$ elements of the sequence $x$, i.
e. on $x_{1},x_{2},\cdots,x_{2k+1}$. Now we can state Theorem 1 of
Turyn's and Storer's paper ``On Binary Sequences'' \cite{turyn1961binary}.
\begin{thm}
Let $x$ satisfy equation (k) for $1\leq k\leq t.$ Let $x_{i}$=1
for \textup{$1\leq i\leq p$, $x_{p+1}=-1$; assume $p>1$. Then}

(i) \hspace{2em} \textup{$x_{i}x_{i+1}=x_{2i}x_{2i+1}$ for $1\leq i\leq t$}

(ii)\hspace{2em} $p\leq2t+1$ implies $p$ is odd

(iii)\hspace{2em} $pj+r\leq2t+1$, \textup{$1\leq r\leq p$} implies
\textup{$x_{p(j-1)+r}=x_{p(j-1)+1}$}

(iv)\hspace{2em} $z_{j}=x_{p(j-1)+1}$ satisfy equations (k) for
$k\leq\frac{t}{p}$
\end{thm}

\section{Counterexamples to Theorem 1 }

In this section we show that there are many counterexamples to Theorem
1 (iv). Consider a binary sequence $x$ with length $n>19$, $x_{1}=1$
and whose first 19 elements have the run length encoding $(3,3,6,3,2,2)$.
For such a sequence $x$ we then have $x_{1}=x_{2}=x_{3}=1$, $x_{4}=x_{5}=x_{6}=-1$,
$x_{7}=x_{8}=\cdots=x_{12}=1$, $x_{13}=x_{14}=x_{15}=-1$, $x_{16}=x_{17}=1$
and $x_{18}=x_{19}=-1$. Let $p$ and $z$ be as in Theorem 1. Hence
$p=3$ and $z=(1,-1,1,1,-1,1,-1,\cdots)$. Furthermore, it is not
difficult to show that in this case $x$ satisfies equation (\emph{k})
for $1\leq k\leq9$; by setting $t=9$ all assumptions of Theorem
1 are therefore satisfied. Now Theorem 1 (iv) claims that the binary
sequence $z$ satisfies equation (\emph{k}) for all $1\leq k\leq\tfrac{t}{p}=3$.
However, this is not true: $z$ does not satisfy equation (\emph{k})
for $k=3$.

As a further example consider a binary sequence $x$ of length $n>33$,
$x_{1}=1$ and whose first 33 elements have the run length encoding
$(5,5,10,5,4,4)$. Again let $p$ and $z$ be as in Theorem 1. Hence
$p=5$ and $z=(1,-1,1,1,-1,1,-1,\cdots)$. Hence $x$ satisfies equation
(\emph{k}) for $1\leq k\leq16$; therefore, we can set $t=16$ in
Theorem 1. Again, Theorem 1 \emph{(iv)} claims that the binary sequence
$z$ satisfies equation (\emph{k}) for all $1\leq k\leq\tfrac{t}{p}=3.2$.
However, this is not true: $z$ does not satisfy equation (\emph{k})
for $k=3$. 

Similar counterexamples for other values of $p\geq3$ can be found
by analyzing sequences whose run length encodings start with $(p,p,2p,p,p-1,p-1,\cdots)$.
Other patterns like $(p,p,p,p,2p,2p,2p-1,p-1,\cdots)$ or $(p,p,2p,p,3p,p,p-1,1,3\cdots)$
also produce counterexamples to Theorem 1. For example, a counterexample
for $p=5$ is given by a binary sequences $x$ of length $n>53$ whose
first 53 elements have the run encoding $(5,5,5,5,10,10,9,4)$ or
$(5,5,10,5,15,5,4,1,3)$; then $x$ satisfies equation (\emph{k})
for all for $1\leq k\leq26$. But although $5\leq\tfrac{26}{p}$,
the corresponding sequence $z$ does not satisfy equation (\emph{k})
for $k=5$ which contradicts Theorem 1 \emph{(iv)}.

\section{What about Theorem 2? }

Theorem 2 in \cite{turyn1961binary} says that there exists no Barker
sequence of odd sequence length $n>13$. It is well-known and often
cited. However, in the proof of Theorem 2 Turyn and Storer use their
Theorem 1 which is (at least in parts) contradicted by the above counterexamples.
An analysis of the proof of Theorem 2 shows that the proof relies
heavily on Theorem 1 \emph{(iii)} but not explicitly on Theorem 1
(\emph{iv}). However, in \cite{turyn1961binary} by induction on $t$
the statement \emph{(iii) }of Theorem 1 is proved simultaneously with
the statement (\emph{iv}) of Theorem 1. Hence, the prove of Theorem
1 \emph{(iii) }also relies on the incorrect Theorem 1 \emph{(iv)};
thus the proof of Theorem 2 in \cite{turyn1961binary} is not correct.
It is not at all clear how Theorem 1 and its proof must be modified
in \cite{turyn1961binary} in order to get a correct proof of the
well-known Theorem 2. An alternative proof of Theorem 2 can for example
be found in \cite{borwein2013note}. Although in many respects quite
similar to the original proof of Turyn and Storer the in \cite{borwein2013note}
presented proof does not rely on the Theorem 1 \emph{(iii)} or Theorem
1 \emph{(iv)}; instead Newton's identities are used in order to prove
a result which is quite similar to Theorem 1 (iii). 

\bibliographystyle{ieeetr}
\bibliography{barker_jw1}

\bigskip{}

\end{document}